\documentclass[12 pt]{amsart}

\usepackage{graphicx}
\usepackage{color}
\usepackage{amssymb, amsmath, amsthm}
\usepackage{mathptmx}

\usepackage{epsfig}
\usepackage{epstopdf}
\usepackage{graphicx}

\newtheorem{thm}{Theorem}[section]

\newtheorem{cor}[thm]{Corollary}

\newtheorem{defn}[thm]{Definition}
\newtheorem{lemma}[thm]{Lemma}

\newtheorem{conj}{Conjecture}

\newcommand{\R}{\mathbb R} 

\newcommand{\bi}{\begin{itemize}}
\newcommand{\ei}{\end{itemize}}
\newcommand{\be}{\begin{enumerate}}
\newcommand{\ee}{\end{enumerate}}

\newcommand{\Aa}{\mathcal{A}}
\newcommand{\B}{\mathcal{B}}

\newcommand{\emp}{\emptyset}

\newcommand{\A}{\alpha}
\newcommand{\pd}{\partial}

\newcommand{\tr}{\text{trunk}}

\begin{document}
\title{A lower bound on the width of satellite knots}
\author{Alexander Zupan}

\maketitle
\begin{abstract}
Thin position for knots in $S^3$ was introduced by Gabai in \cite{Gabai} and has been used in a variety of contexts.  We conjecture an analogue to a theorem of Schubert and Schultens concerning the bridge number of satellite knots.  For a satellite knot $K$, we use the companion torus $T$ to provide a lower bound for $w(K)$, proving the conjecture for $K$ with a 2-bridge companion.  As a corollary, we find thin position for any satellite knot with a braid pattern and 2-bridge companion.
\end{abstract}

\section{Introduction}

Thin position for knots in $S^3$ was introduced by Gabai in \cite{Gabai} and has since been studied extensively.  Although thin position has been used in a variety of different proofs, there are relatively few methods for putting specific knots into thin position.  Thin position of a knot always provides a useful surface; either a level sphere is a bridge sphere for the knot or the thinnest thin sphere is incompressible in the complement of the knot, as shown by Wu \cite{wu}. \\

In some sense, width can be considered to be a refinement of bridge number, although recently it has been shown in \cite{addit} that one can not recover the bridge number of a knot $K$ from the thin position of $K$.  On the other hand, if $K$ is small, then $w(K) = 2 \cdot b(K)^2$ and any thin position of $K$ is a bridge position.  In his classic paper on the subject \cite{schub}, Schubert proved that for any two knots $K_1$ and $K_2$, $b(K_1 \# K_2) = b(K_1) + b(K_2) - 1$.  This was later reproved by Schultens in \cite{bridge}. \\

Unfortunately, we cannot hope for a similar statement to hold for width.  In \cite{planar}, Scharlemann and Schultens establish $\max \{w(K_1),w(K_2)\}$ as a lower bound for $w(K_1 \# K_2)$, and Blair and Tomova prove that this bound is tight in some cases \cite{addit}, while Rieck and Sedwick \cite{stack} demonstrate that the bound is not tight for small knots.  Both Schubert and Schultens also prove the following:
\begin{thm}
Let $K$ be a satellite knot with pattern $\hat{K}$ and companion $J$, where $n$ is the winding number of $\hat{K}$.  Then
\[ b(K) \geq n \cdot b(J).\]
\end{thm}
We make an analogous conjecture:
\begin{conj}
Let $K$ be a satellite knot with pattern $\hat{K}$ and companion $J$, where $n$ is the winding number of $\hat{K}$.  Then
\[ w(K) \geq n ^2\cdot w(J).\]
\end{conj}
In this paper, we provide a weaker lower bound for $w(K)$.  Our main theorem is as follows:
\begin{thm}
Let $K$ be a satellite knot with pattern $\hat{K}$, where $n$ is the winding number of $\hat{K}$.  Then
\[ w(K) \geq 8n^2.\]
\end{thm}
This proves the conjecture in the case that the companion $J$ is a 2-bridge knot, since the width of such $J$ is 8.  As a corollary, if $K$ is a satellite with a 2-bridge companion and its pattern $\hat{K}$ is a braid with index $n$, then any thin position is a bridge position for $K$.

\section{Preliminaries}
Let $K$ be a knot in $S^3$, and let $\mathcal{M}(K)$ denote the collection of Morse functions $h: S^3 \rightarrow \R$ with exactly two critical points on $S^3$, denoted $\pm \infty$, and such that $h \mid_K$ is also Morse.  (Equivalently, we could fix some Morse function $h$ and look instead at the collection of embeddings of $K$ into $S^3$.)  For every $h \in \mathcal{M}(K)$, let $c_0 < c_1 < \dots < c_n$ denote the critical values of $h \mid_K$.  Choose regular levels $c_0 < r_1 < c_1 < \dots < r_n < c_n$, and define
\begin{eqnarray*}
w(h) &=& \sum_{i=1}^n |K \cap h^{-1}(r_i)|, \\
b(h) &=& \frac{n+1}{2}, \\
\tr(h) &=& \max |K \cap h^{-1}(r_i)|.
\end{eqnarray*}
Now, let
\begin{eqnarray*}
w(K) &=& \min_{h \in \mathcal{M}(K)} w(h),\\
b(K) &=& \min_{h \in \mathcal{M}(K)} b(h), \\
\tr(K) &=& \min_{h \in \mathcal{M}(K)} \tr(h).
\end{eqnarray*}
These three knot invariants are called the width, bridge number, and the trunk of $K$, respectively.  Width was defined by Gabai in \cite{Gabai}, and trunk was defined by Ozawa in \cite{waist}.  Observe that $b(K)$ is the least number of maxima of any embedding of $K$.  If $h \in \mathcal{M}(K)$ satisfies $w(K) = w(h)$, we say that $h$ is a thin position for $K$.  If $h \in \mathcal{M}(K)$ satisfies $b(K) = b(h)$ and all maxima of $h\mid_K$ occur above all minima, then we say that $h$ is a bridge position for $K$. \\

In \cite{planar}, the authors give an alternative formula for computing width, which involves thin and thick levels.  Let $h \in \mathcal{M}(K)$ with critical and regular values as defined above.  Then we say $h^{-1}(r_i)$ is a thick level if $|K \cap h^{-1}(r_i)| > |K \cap h^{-1}(r_{i-1})|,|K \cap h^{-1}(r_{i+1})|$ and $h^{-1}(r_i)$ is a thin level if $|K \cap h^{-1}(r_i)| < |K \cap h^{-1}(r_{i-1})|,|K \cap h^{-1}(r_{i+1})|$, where $1 < i < n$.  Note that if $h$ is a bridge position for $K$, then $h$ has exactly one thick level and no thin levels.  Letting $a_1,\dots,a_m$ denote the number of intersections of the thick levels with $K$ and $b_1,\dots,b_{m-1}$ denote the number of intersection of the thin levels with $K$, the width of $h$ is given by
\[ w(h) = \frac{1}{2} \left( \sum_{i=1}^m a_i^2 - \sum_{i=1}^{m-1} b_i^2\right).\]
In particular, we see that for every $h \in \mathcal{M}(K)$, there exists $a_i \geq \tr(K)$, which implies that
\[ w(K) \geq \frac{\tr(K)^2}{2}.\]

The knots we will be concerned with are satellite knots, defined below:
\begin{defn}
Let $\hat{K} \subset V$ be a knot contained in a solid torus $V$ with core $C$ and such that every meridian of $V$ intersects $\hat{K}$, and let $J$ be any nontrivial knot.  Suppose that $\varphi:V \rightarrow S^3$ is an embedding such that $\varphi(C)$ is isotopic to $J$ in $S^3$.  Then $K = \varphi(\hat{K})$ is called a satellite knot with companion $J$ and pattern $\hat{K}$.
\end{defn}
Essentially, to construct a satellite knot $K$, we start with a pattern in a solid torus and then tie the solid torus in the shape of the companion $J$.  We will need several more definitions to state the main result:
\begin{defn}
Let $\hat{K}$ be a pattern contained in a solid torus $V$.  The winding number of $\hat{K}$, $\#(\hat{K})$, is the absolute value of the algebraic intersection number of any meridian disk of $V$ with $\hat{K}$.
\end{defn}
Equivalently, if $\A:S^1 \rightarrow V$ is an embedding such that $\A(S^1) = \hat{K}$ and $r:V \rightarrow S^1$ is a strong deformation retract of $V$ onto its core, then $\#(\hat{K})$ agrees with the degree of the map $r \circ \A$.
\begin{defn}
Let $\hat{K}$ be a pattern contained in a solid torus $V$.  We say that $\hat{K}$ is a braid of index $n$ if there is a foliation of $V$ such that every leaf is a meridian disk intersecting $\hat{K}$ exactly $n$ times.
\end{defn}
In the case that $\hat{K}$ is a braid of index $n$, it is clear that $\#(\hat{K}) = n$.  For an example, consider Figure 1.  On the left, we see a braid pattern of index $3$, $\hat{K}$, contained in a solid torus $V$.  On the right, $V$ is embedded in such a way that its core is a trefoil.  Thus, the knot $K$ on the right is a satellite knot with trefoil companion and pattern $\hat{K}$.
\begin{figure}
  \centering
      \includegraphics[width=1.0\textwidth]{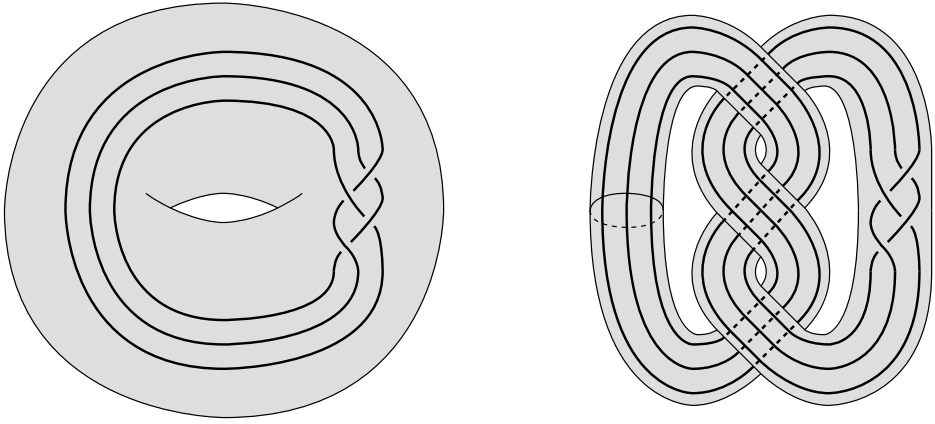}
  \caption{On the left, pattern $\hat{K}$ is shown contained in a solid torus.  On the right, we see a satellite knot $K$ with pattern $\hat{K}$ and trefoil companion.}
\end{figure}

\section{Reducing the saddle points on the companion torus}
From this point on, we set the convention that $K$ is a satellite knot with companion $J$ and pattern $\hat{K}$ contained in a solid torus $\hat{V}$, $\varphi$ is an embedding of $\hat{V}$ into $S^3$ that takes a core of $\hat{V}$ to $K'$, $V= \varphi(\hat{V})$, and $T = \pd V$.  Further, we will let $h \in \mathcal{M}(K)$ and perturb $V$ slightly so that $h \mid_T$ is Morse.  We wish to restrict our investigation to tori $T$ with only certain types of saddle points.  In this vein, we follow \cite{bridge}, from which the next definition is taken.
\begin{defn}
Consider the singular foliation, $F_T$, of $T$ induced by $h \mid_T$.  Let $\sigma$ be a leaf corresponding to a saddle point.  Then one component of $\sigma$ is the wedge of two circles $s_1$ and $s_2$.  If either is inessential in $T$, we say that $\sigma$ is an inessential saddle.  Otherwise, $\sigma$ is an essential saddle.
\end{defn}
The next lemma is the Pop Over Lemma from \cite{bridge}:
\begin{lemma}
If $F_T$ contains inessential saddles, then after a small isotopy of $T$, there is an inessential saddle $\sigma$ in $T$ such that
{\be
\item $s_1$ bounds a disk $D_1 \subset T$ such that $F_T$ restricted to $D_1$ contains only one maximum or minimum,
\item for $L$ the level surface of $h$ containing $\sigma$, $D_1$ co-bounds a 3-ball $B$ with a disk $\tilde{D}_1 \subset L$ such that $B$ does not contain $\pm \infty$ and such that $s_2$ lies outside of $\tilde{D}_1$.
\ee}
\end{lemma}
In the following lemma, we mimic Lemma 2 of \cite{bridge} with a slight modification to preserve the height function $h$ on $K$:
\begin{lemma}
There exists an isotopy $f_t:S^3 \rightarrow S^3$ such that $f_0 = \text{id}$, $h = h \circ f_1$ on $K$, and the foliation of $T$ induced by $h \circ f_1$ contains no inessential saddles.
\begin{proof}
Suppose that $T$ has an inessential saddle, $\sigma$, lying in the level 2-sphere $L$.  By the previous lemma, we may suppose that $\sigma$ is as described above, and suppose without loss of generality that $D_1$ contains only one maximum.  By slightly pushing $D_1$ into $\text{int}(B)$, we can create a new closed ball $B'$ such that $B' \cap D_1 = \emp$ and $(K \cup T) \cap \text{int}(B) \subset B'$.  First, we isotope $B'$ vertically until it lies below $L$, and then isotope $D_1$ down until the maximum of $D_1$ cancels out the saddle point $\sigma$.  Now, there exists a monotone increasing arc beginning at the highest point of $B'$, passing through the disk $\tilde{D}_2$ bounded by $s_2$, intersecting only maxima of $T$, and disjoint from $K$.  Thus, we may isotope $B'$ vertically through a regular neighborhood of $\A$, increasing the heights of maxima of $T$ if necessary, until the heights of maxima and minima of $K \cap \text{int}(B')$ are the same as before any of the above isotopies.  We see that after isotopy $T$ has one fewer inessential saddle and no new critical points have been created.  See Figure 2.  Repeating this process, we eliminate all inessential saddles via isotopy.
\end{proof}
\end{lemma}
\begin{figure}
  \centering
      \includegraphics[width=1.1\textwidth]{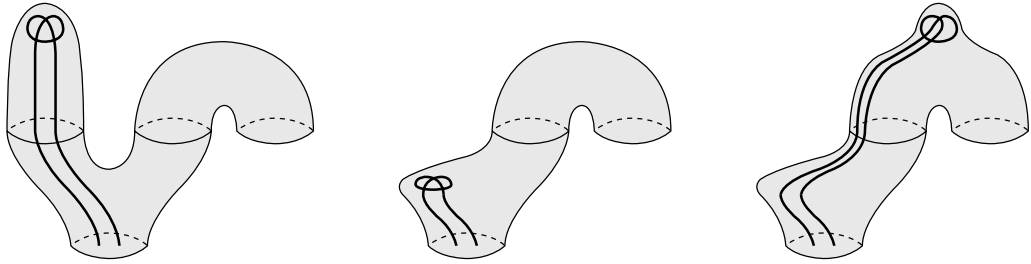}
  \caption{First, we cancel the inessential saddle, shown center.  Then we isotope any part of $K$ or $T$ contained in $B$ along an increasing arc $\A$, increasing maxima of $T$ is necessary, so that $h\mid_K$ is unchanged with respect to the end product of our isotopy, shown at right.}
\end{figure}
Thus, from this point forward, we may replace any $h \in \mathcal{M}(K)$ with $h \circ f_1$ from the lemma without changing the information carried by $h\mid_K$; thus we may suppose that the torus $T = \partial V$ contains no inessential saddles.  It follows that if $\gamma$ is a loop contained in a level 2-sphere that bounds a disk $D \subset T$, then $D$ contains exactly one critical point, a minimum or a maximum.  If not, $D$ must contain a saddle point, which is necessarily inessential.

\section{The connectivity graph}
For each regular value $r$ of $h\mid_{T,K}$, we have that $h^{-1}(r)$ is a level 2-sphere $S^2$ and $h^{-1}(r) \cap T$ is a collection of simple closed curves.  Let $\gamma_1,\dots,\gamma_n$ denote these curves. \\

A bipartite graph is a graph together with a partition of its vertices into two sets $\Aa$ and $\B$ such that no two vertices from the same set share an edge.  We will create a bipartite graph $\Gamma_r$ from $h^{-1}(r)$ as follows: Cut the 2-sphere $h^{-1}(r)$ along $\gamma_1,\dots,\gamma_n$, splitting $h^{-1}(r)$ into a collection of planar regions $R_1,\dots,R_m$.  The vertex set $\{v_1,\dots,v_m\}$ of $\Gamma_r$ corresponds to the regions $R_1,\dots,R_m$, and the edges correspond to the curves $\gamma_1,\dots,\gamma_n$ that do not bound disks in $T$.  For each such $\gamma_i$, make an edge between $v_j$ and $v_k$ if $\gamma_i = R_j \cap R_k$ in $h^{-1}(r)$.  To see that $\Gamma_r$ is bipartite, we create two vertex sets $\Aa_r$ and $\B_r$, letting $v_i \in \Aa_r$ if $R_i \subset V$, and $v_i \in \B_r$ otherwise.  We call $\Gamma_r$ the \textbf{essential connectivity graph} with respect to the regular value $r$ of $h$, where the word ``essential" emphasizes the fact that edges correspond to only those $\gamma_i$ that are essential in $T$.  Note that since each $\gamma_i$ separates $h^{-1}(r)$, the graph $\Gamma_r$ must be a tree.  An endpoint of $\Gamma_r$ is a vertex that is incident to exactly one edge. \\

\begin{figure}
  \centering
      \includegraphics[width=1.0\textwidth]{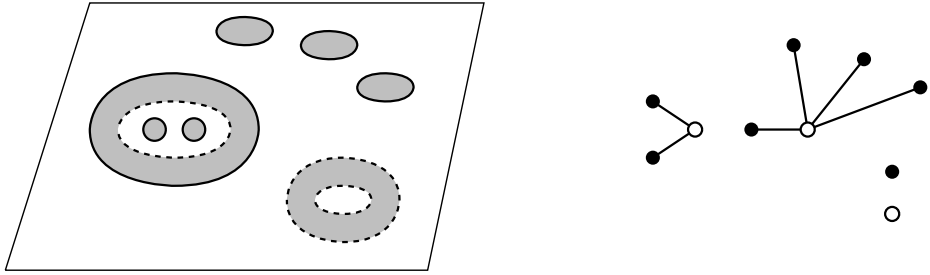}
  \caption{We see a level 2-sphere at left and its corresponding essential connectivity graph at right.  Note that dotted curves on the left correspond to curves bounding disks in $T$.}
\end{figure}
For instance, in Figure 3 we see a possible level 2-sphere and corresponding essential connectivity graph.  Observe that since $V$ is a knotted solid torus, $T$ is only compressible on one side, and every compression disk for $T$ is a meridian of $V$.  This leads to the third lemma:

\begin{lemma}
If $v_i \in \Gamma_r$ is an endpoint, then $v_i \in \Aa_r$.
\begin{proof}
Suppose $R_i$ is the region in $h^{-1}(r)$ corresponding to $v_i$.  Then $\pd R_i$ contains exactly one essential curve in $T$, call it $\gamma$, and some (possibly empty) set of curves that bound disks in $T$.  Since each of these disks contains only one maximum or minimum by the discussion above, any two must be pairwise disjoint.  Thus, we can glue each disk to $R_i$ to create an embedded disk $D$ such that $\partial D = \gamma$.  Now, push each glued disk into a collar of $T$ in $V$, so that $T \cap \text{int}(D) = \emp$, and thus $D$ is a compression disk for $T$.  We conclude $D \subset V$ and $R_i \cap D \neq \emp$, implying $R_i \subset V$ and $v_i \in \Aa_r$.
\end{proof}
\end{lemma}
Using similar arguments, we prove the next lemma:
\begin{lemma}
Suppose that $v_1,\dots,v_n \subset \Gamma_r$ are endpoints corresponding to regions $R_1,\dots,R_n \subset h^{-1}(r)$, where each $R_i$ contains exactly one curve $\gamma_i$ that is essential in $T$.  Then $\gamma_1,\dots,\gamma_n$ bound meridian disks $D_1,\dots,D_n \subset V$ such that $K \cap D_i \subset R_i$ for all $i$.
\begin{proof}
The existence of the disks $D_1,\dots,D_n$ is given in the proof of Lemma 3.  Thus, suppose that $\Delta$ is a disk glued to $R_i$ to construct $D_i$.  When we push $\Delta$ into a collar of $T$, we can choose this collar to be small enough so that it does not intersect $K$.  Thus, we may suppose that $\Delta \cap K = \emp$ for every such $\Delta$, which implies that all intersections of $K$ with $D_i$ must be contained in $R_i$.
\end{proof}
\end{lemma}
We note that the Lemmas 3 and 4 are inspired by the proof of Theorem 1.9 of \cite{waist}.  Essentially, Lemma 4 demonstrates that even though the set of meridian disks $D_1,\dots,D_n$ may not be level, we may assume they are level for the purpose of counting intersections of $K$ with $h^{-1}(r)$, since any intersection of $K$ with one of these disks occurs in one of the level regions $R_i$.  Hence, we define the trunk of a level 2-sphere.
\begin{defn}
Let $r$ be a regular value of $h\mid_{T,K}$.  We define the trunk of the level 2-sphere $h^{-1}(r)$, denoted $\tr(r)$, to be the number of endpoints of $\Gamma_r$.
\end{defn}
For example, if $r$ is the regular value whose essential connectivity graph is pictured in Figure 3, then $\tr(r) = 6$.  We are now in a position to use the winding number of the pattern $\hat{K}$.
\begin{lemma}
Let $r$ be a regular value of $h \mid_{T,K}$.
{\bi
\item If $\tr(r)$ is even, then $|K \cap h^{-1}(r)| \geq \#(\hat{K}) \cdot \tr(r)$;
\item if $\tr(r)$ is odd, then $|K \cap h^{-1}(r)| \geq \#(\hat{K}) \cdot [\tr(r) + 1]$.
\ei}
\begin{proof}
First, suppose that $m = \tr(r)$ is even and let $n = \#(\hat{K})$.  Since each meridian of $V$ has algebraic intersection $\pm n$ with $K$, we know that each meridian must intersect $V$ in at least $n$ points.  Let $v_1,\dots,v_m$ be endpoints of $\Gamma_r$ corresponding to regions $R_1,\dots,R_m$.  By Lemma 4, $|K \cap R_i| = |K \cap D_i| \geq n$ for each $i$.  Further, since these regions are pairwise disjoint, it follows that $|K \cap h^{-1}(r)| \geq n \cdot m$, completing the first part of the proof. \\

Now, suppose that $m$ is odd.  If $N_1$ is the algebraic intersection number of $K$ with $R = \cup R_i$, we have that
\[ N_1 = \sum_{i=1}^m \pm n.\]
In particular, as $m$ is odd it follows that $|N_1| \geq n$.  Let $R' = \overline{h^{-1}(r) - R}$.  Then $R' \cap R \subset T$, so $K$ does not intersect $R' \cap R$.  Let $N_2$ denote the algebraic intersection number of $K$ with $R'$.  Since $h^{-1}(r)$ is a 2-sphere which bounds a ball in $S^3$, $h^{-1}(r)$ is homologically trivial, implying that the algebraic intersection of $K$ with $h^{-1}(r)$ is zero.  This means $N_1 + N_2 = 0$, so $|N_2| \geq n$ and thus $|K \cap R'| \geq n$.  Finally, putting everything together, we have
\[ |K \cap h^{-1}(r)| = |K \cap R| + |K \cap R'| = \sum_{i=1}^m |K \cap R_i| + |K \cap R'| \geq n \cdot (m+1).\]
\end{proof}
\end{lemma}

\section{Bounding the width of satellite knots}
We will use the trunk of the level surfaces to impose a lower bound on the trunk of a $K$, which in turn forces a lower bound on the width of the $K$.  We need the following lemma, which is Claim 2.4 in \cite{waist}:
\begin{lemma}
Let $S$ be a torus embedded in $S^3$, and let $h:S^3 \rightarrow \R$ be a Morse function with two critical points on $S^3$ such that $h \mid_S$ is also Morse.  Suppose that for every regular value $r$ of $h \mid_S$, all curves in $h^{-1}(r)\cap S$ that are essential in $S$ are mutually parallel in $h^{-1}(r)$.  Then $S$ bounds solid tori $V_1$ and $V_2$ in $S^3$ such that $V_1 \cap V_2 = T$.
\end{lemma}
As a result of this lemma, we have
\begin{cor}
There exists a regular value $r$ of $h \mid_{T,K}$ such that $\tr(r) \geq 3$.
\begin{proof}
Suppose not, and let $r$ be any regular value of $h \mid_{T,K}$ such that $h^{-1}(r)$ contain essential curves in $T$.  Such a regular value must exist; otherwise $T$ could not contain a saddle point.  By assumption, $\tr(r) \leq 2$, so $\Gamma_r$ has exactly two endpoints, $v_1$ and $v_2$.  But this implies that $\Gamma_r$ is a path, and thus all essential curves in $h^{-1}(r)$ are mutually parallel.  As this is true for every such regular value $r$, we conclude by Lemma 6 that $V$ is an unknotted solid torus, contradicting the fact that $K$ is a satellite knot with nontrivial companion $J$.
\end{proof}
\end{cor}
This brings us to our main theorem.
\begin{thm}
Suppose $K$ is a satellite knot with pattern $\hat{K}$, where $n = \#(\hat{K})$.  Then
\[ w(K) \geq 8n^2.\]
\begin{proof}
Choose a height function $h \in \mathcal{M}(K)$ such that $\tr(h) = \tr(K)$.  Since $K$ is a satellite knot, $K$ is contained in a knotted solid torus $V$.  Let $T = \pd V$, and if necessary perturb $T$ slightly so that $h \mid_T$ is also Morse.  By Corollary 1 above, there exists a regular value $r$ of $h$ such that $\tr(r) \geq 3$.  From Lemma 5, it follows that $|K \cap h^{-1}(r)| \geq 4n$.  Since $\tr(K) = \tr(h)$, and $\tr(h)$ corresponds to the level of $h$ with the greatest number of intersections with $K$, we have $\tr(h) \geq 4n$.  Finally, using the lower bound for width based on trunk,
\[ w(K) \geq \frac{\tr(K)^2}{2} \geq 8n^2,\]
as desired.
\end{proof}
\end{thm}

\begin{cor}
Suppose $K$ is a satellite knot, with pattern $\hat{K}$ and companion $J$.  If $\hat{K}$ is a braid of index $n$ and $J$ is a 2-bridge knot, then $w(K) = 8n^2$ and any thin position for $K$ is a bridge position.
\begin{proof}
For such $K$ we can exhibit a Morse function $h \in \mathcal{M}(K)$ such that $w(h) = 8n^2$, $b(h) = 2n$, and $\tr(h) = 4n$.  By \cite{bridge}, $b(K) = b(h)$, so $h$ is both a bridge and thin position for $h$, and further every bridge position $h'$ for $K$ satisfies $w(h') = 8n^2$ and is also thin.  It follows from the proof of the above theorem that $\tr(K) = 4n$, so any $h \in \mathcal{M}(K)$ that is not a bridge position satisfies $w(h) > 8n^2$.
\end{proof}
\end{cor}
\section{Acknowledgements}
I would like to thank Maggy Tomova and Charlie Frohman for their support and advice, including numerous helpful conversations.


\begin{thebibliography}{3}
\bibitem{addit} R Blair, M Tomova, Width is not additive, preprint, available at http://arxiv.org/abs/1005.1359.
\bibitem{Gabai} D Gabai, Foliations and the topology of 3-manifolds III, J. Differential Geom. 26 
(1987) 479-536.
\bibitem{waist} M Ozawa, Waist and trunk of knots, Geom. Ded. (2010) Online First.
\bibitem{stack} Y Rieck, E Sedwick, Thin position for a connected sum of small knots, Alg. Geom. Top. 2 (2002) 297-309.
\bibitem{planar} M Scharlemann, J Schultens, 3-manifolds with planar presentations and the width of satellite knots, Trans. Amer. Math. Soc. 358 (2006) 3781-3805.
\bibitem{schub} H Schubert, \"{U}ber eine numerische Knoteninvariante, Math. Z. 61 (1954) 245-288.
\bibitem{bridge} J Schultens, Additivity of bridge numbers of knots, Math. Proc. Cambridge Philos. Soc. 135 (2003) 539-544.
\bibitem{wu} Y-Q Wu, Thin position and essential planar surfaces, Proc. Amer. Math. Soc. 132 (2004), no. 11, 3417-3421 (electronic).
\end{thebibliography}
\end{document}